\documentclass[11pt]{article}
\usepackage{amssymb}
\usepackage{latexsym}

\newtheorem{thm}{Theorem}[section]
\newtheorem{prop}[thm]{Proposition}
\newtheorem{lemma}[thm]{Lemma}

\begin{document}

\title{\bf  Symmetries of 3-webs around a point}

\author{ Jean Paul Dufour. \\}

\maketitle

\begin{abstract}
Let $W$ be a real analytic planar 3-web defined on a neighborhood of a point $M$. In this paper whe consider local diffeomorphisms of the plane which preserve $M$ and which map any foliation of $W$ onto a (not necessarily the same) foliation. We call them "symmetry" of $W$ around $M$. There is three types of such symmetries: the first one are those which preserve each foliation; the second one are those which preserve only one foliation and permute the two others; the third are those which permute circularly the three foliations. 

We know that hexagonal (i.e. flat) 3-webs have always the three types of symmetries around each point. In this paper we study the non-flat case.

 We give a classification of 3-webs which admit a symmetry of first or second type. The problem is more difficult for the third case but we give a method for building all 3-webs  with this type of symmetry. We also give a precise example which is non-flat.

\end{abstract}


\section{Introduction}

In this text we work in the real analytic framework. The results could be extended to the smooth and to the complex frameworks.

A planar 3-web is a triple  of 1-dimensional foliations, mutually tranversal, on an open domain of the plane. In the sequel we will work only with these 3-webs: so we forget the word {\sl planar}. Let $W$ be such a web and $M$ a point of  its domain.

A {\bf symmetry of $W$ around $M$} is a local diffeomorphism of the plane which preserves $M$ and maps any foliation of $W$ to another (possibly different) foliation.

We have to consider three  different types  of such symmetries:    {\bf simple symmetries}  are those which preserve each foliation;  {\bf mirror symmetries}  are those which preserve only one foliation and permute the two others;  {\bf circular symmetries} are those which give a circular permutation of the three foliations.

We choose local coordinates $(x,y)$ near $M$, vanishing at this point, such that the foliations  are given by the verticals lines $ x={ constant},$ the horizontals lines $y={ constant}$ and  the level sets of some function $f$ (the sets of $(x,y)$ such that $ f(x,y)
={constant}).$ Such a 3-web is denoted $(x,y,f).$

We know that $W$ is {\bf flat} near $M$ iff we can find such coordinates  with $f(x,y)=x+y.$ In that case it is easy to see that it has  circular, mirror or non-trivial simple symmetries.

\section{Normal forms and simple symmetries.}

A classical result (\cite{DJ}), based on a result of S. Sternberg (\cite{SS}),
  claims that we can always choose coordinates $(x,y)$ vanishing at $M$ such that
\begin{equation} f(x,y)=x+y +xy(x-y)g(x,y)\end{equation}
where $g$ is an analytic (or smooth) function. 

We say that $ (x,y, x+y +xy(x-y)g(x,y))$ is a  {\bf normal form } of $W$ near $M$. They are characterised by the three equations 
$$f(t,0)=t,\\\ f(0,t)=t,\\\ f(t,t)=2t,$$ for any small $t,$ and
 $W$ is {\bf non-flat} (\cite{AH}) near $M$ iff the function $g(x,y)$ does not vanish
 identically near the origin.
Moreover these coordinates are unique up to a homothety $(x,y)\mapsto (\lambda x,\lambda y)$. At least this notion proves easily the following result.

\begin{prop} If $W$ admits a simple symmetry $\phi$  around $M$ with has $Id$ (identity) as tangent map at $M$ then $\phi$ is the identity.\end{prop}

We denote by $\sum g_{rs}x^ry^s$ the  Taylor expansion of the fonction $g$ near the origin, the above homothety maps $(x,y,f(x,y))$ to a  the web $(x,y,x+y+xy(x-y)G(x,y))$ with 
$$G(x,y)=\sum g_{rs}\lambda^{r+s+2}x^ry^s.$$

Suppose there is at least one non-vanishing coefficient $ g_{r_0s_0}$
  with $r_0+s_0$ odd. 
When we choose correctly $\lambda$, we obtain that the new coefficient of $x^{r_0}y^{s_0}$ is $1.$ With this new constraint, the normal form is unique and the only simple symmetry is the identity ($Id$).

Suppose now that all the coefficients $g_{rs}$ with $r+s$ odd vanish. Then the homothety $(x,y)\mapsto (-x,-y)$ ($-Id$) doesn't change our normal form. So we can conclude with the following 
\begin{thm} If the 3-web $W$ is non-flat, either it has only $Id$ as simple symmetry around $M$, or it has two simple symmetries $Id$ and  $\phi$, where $\phi$ is a non trivial involution ($\phi^2=I$ and the linear part of $\phi$ at the origin is $-Id$).

Moreover it has exactly two simple symmetries if and only if its normal form  $(x,y, x+y +xy(x-y)g(x,y))$ is such that $g(x,y)=g(-x,-y)$ for every $x$ and $y.$\end{thm}

As far I know, this result never appeared in the classical web litterature.

 {\bf Remark:} The notion of normal form works also in the complex case. We have an analog of the preceeding proposition but we may have more than the above non-trivial simple symmetries: We suppose that $n$ is the lowest integer such that they are non-zero $g_{rs}$ with $n=r+s$.  For a convenient choice of $\lambda$ the diffeomorphism $\lambda . Id$ gives a new normal form but with $g_{r_0s_0}=1$ for some couple $(r_0,s_0)$ with $r_0+s_0=n.$ The only diffeomorphisms $\phi$ which preserve this constraint are of the form $\mu . Id$ where $\mu$ is a $(n+2)$-root of the unit ($\mu^{n+2}=1$).
As $\phi$ preserves only the monomial $x^ry^s$  with $r+s+2=k.(n+2)$ for any  integer $k$, $\mu.Id$ will give $n+1$  non-trivial simple symmetries when $G$ contains only these monomials in its Taylor expansion.

\section{Mirror symmetries of 3-webs.}

     Suppose that $W$ has a mirror symmetry $\phi$ around $M.$ 

We  first remark that $\phi^2$ is a simple symmetry; the theorem above implies that we have two possibilities : either $\phi^2=Id$ or $\phi^4=Id$ and the linear part of $\phi^2$ at the origin is $-Id.$ 

We can choose local coordinates $(x,y)$, which give the normal form for $W$, and such that $\phi$ permutes the foliations $ x={ constant}$ and  $y={ constant}$. Thus $\phi$ has the form $(x,y)\mapsto (X(y),Y(x))$ for two diffeomorphisms $X$ and $Y$  of the line (fixing the origin). The remark above proves that there are two cases: the linear part of $\phi$ is $(x,y)\mapsto (y,x)$ or  $(x,y)\mapsto (-y,-x)$.

 Since $\phi$ is a mirror symmetry, it must preserve the third foliation. As this last is given by the the level lines of a fonction
$f(x,y)=x+y+xy(x-y)g(x,y),$
this may be written as
$$f((X(y),Y(x))=Z(f(x,y)),$$
for some  third diffeomorphism $Z$ of the line.

When we use the relation $f(t,0)=f(0,t)=t$ we obtain 
$X=Y=Z.$

We have also the obvious property $f(t,t)=2t$. So when we write $x=y=t$ in the above formula we obtain
$$X(2t)=2X(t),$$
which implies that $X$, $Y$ and $Z$ are the same homothety: $Id$ or$-Id.$ So we can rewrite our equation as
$$f( y,x)=f(x,y)$$ or $$f(-y,-x)=-f(x,y)$$
This can be written as
$$g(y, x)=-g(x,y).$$ or $$g(-y,-x)=-g(x,y).$$

\begin{thm} If the 3-web $W$ is non flat, it  has zero, one or two mirror symmetries around $M.$

It has a mirror symmetry $\phi$ which is an involution ($\phi^2=Id$) iff it has a normal form $(x,y,x+y+xy(x-y)g(x,y)),$ with a function $g$ which satisfies $g(y,x)=-g(x,y)$. It has a mirror symmetry $\phi$ which $\phi^4=Id$ iff it has a normal form with a function $g$ which satisfies $g(y,x)=-g(-x,-y)$. It may have simultaneously the two types of mirror symmetries.
\end{thm}

\section{Circular symmetries.}

Suppose that $W$ has a circular symmetry $\phi$ around $M.$ We  choose local coordinates $(x,y)$ which provide to $W$ the  form $ (x,y, f(x,y)),$  with $$f(x,y)=x+y +hot,$$ ("$hot$" means higher order terms)  and such that $\phi$
maps  the foliation by vertical lines to the foliation by horizontal lines, the foliation by horizontal lines to the foliation by the level sets of $f$ and the foliation by the level sets of $f$ to the foliation by the vertical lines. 

Remark that $\phi^3$ is a simple symmetry of $W.$ So, if $W$ is non-flat, we have  $\phi^3=Id$ or $\phi^6=Id$.

In the first case the linear part $\phi^{(1)}$ of $\phi$ at the origin has the form $(x,y)\mapsto (-x-y,x)$ and  $(x,y)\mapsto (x+y,-x)$ in the second case.

In the sequel we study only the first case $\phi^3=Id$; the second can be treated exactly as the first one.

In that  case  $\phi$ gives a local analytic  action of $\mathbb Z /3\mathbb Z$ wich fixes the origin. A classical result claims that such local action is conjugate to its linear part. This means that there is a local diffeomorphism $\psi$ (with identity as linear part at the origin) which fixes the origin and such that
$$\phi =\psi^{-1}\circ \phi^{(1)}\circ \psi$$
on a neighborhood of the origin. 
 With notations $\psi (x,y)=(A(x,y),B(x,y))$ and $\psi^{-1} (x,y)=(U(x,y),V(x,y))$ and $\phi (x,y)=(F(x,y),G(x,y))$ we have the relations
\begin{equation}U(A(x,y),B(x,y))-x =0,\end{equation} \label{Eq} 
\begin{equation}V(A(x,y),B(x,y))-y=0, \end{equation} \label{Eq} 
\begin{equation}F(x,y)=U(-A(x,y)-B(x,y),A(x,y)), \end{equation} \label{Eq} 
\begin{equation}G(x,y)=V(-A(x,y)-B(x,y),A(x,y)). \end{equation} \label{Eq}
The linear parts of $A$ and $U$ at the origin is $x$ and the linear parts of $B$ and $V$ at the origin is $y.$

Now we remark that $\phi$ maps vertical lines to horizontal liness if and only if $G(x,y)$ is independant of $y$. So we have the following equation

\begin{equation}V(-A(x,y)-B(x,y),A(x,y))-\mu(x)=0 \end{equation} \label{Eq}

for some function $\mu =x+hot.$

When the functions $\mu$ and $V(x,y)=y+hot$ are given, the theorem of implicit functions claims that the pair of equations (6) and (3) have a suitable local solution $(A,B)$. Then the equation (2) has a solution $U$ and we obtain $F$ by  (4). 

It is obvious that $\phi$ maps  level sets of the function $F$ to vertical lines. Because $\phi^3=Id$, it also obvious that $\phi$ is a circular symmetry of the 3-web $(x,y,F).$ As we have $F=-x-y+hot$ we obtain $f=-F.$

This gives a way to construct all the non-flat webs which admit a circular symmetry. However this method, starting with differents $(\mu,V)$, could give isomorphic 3-webs; it can also give flat webs. In the sequel we will give a $V$ such that, with $\mu=Id$, we obtain a non-flat 3-web in normal form.

\begin{lemma}  We assume that $P(x,y)$ is a function invariant by $\phi^{(1)}$, i.e. such that $P(-x-y,x)=P(x,y),$ with a vanishing  1-jet at the origin. We assume also that $\theta(x,y)$ is the local solution of the equation $$\theta +P(x+\theta,y+\theta)=0.$$  
For $V=y+P$ and $\mu=Id$, we obtain $f(x,y)=-x-y-3\theta$. \end{lemma}

\noindent {\bf Proof:} We write $A=x+\alpha(x,y),$  $B=y+\beta(x,y)$ and $U=x+\gamma(x,y)$ where the 1-jets  at the origin of $\alpha$, $\beta$ and $\gamma$ vanish. Equations (6) and (3) give $B+P(A,B)=y$ and $A+P(A,B)=x$; so we get $\beta=\alpha$ and $\alpha +P(x+\alpha,y+\alpha)=0.$
So we have $\alpha=\theta.$

Now equation (2) gives $ \gamma=P$ and $F=-A-B+P(A,B)$ leads to  $f(x,y)=-x-y-3\theta$, which achieves our proof.

Consider the polynomial $P_0=xy(x-y)(x+y)(2x+y)(x+2y)$. It is  invariant by $\phi^{(1)}$.
\begin{lemma} The equation $\theta+P_0(x+\theta,y+\theta)=0$ have a solution of the shape $\theta=P_0(-1+hot)$.\end{lemma}

To prove this lemma we impose $\theta=P_0(-1+\eta)$ and we apply the theorem of implicit functions to our equation (considered as an equation of the shape $K(x,y,\eta)=0)$.

Now we find that $F$ has the shape $-x-y-3P_0(1+\eta)$. So our web  has the normal form $(x,y,f)$ with $f=-F$ and it is non-flat.

Remark that the level set $f=0$ is the line $y=-x$. The level set $f=t$ contains the points $(-t,2t),$ $(0,t),$ $(t/2,t/2),$ $(t,0)$ and $(2t,-t)$; they are all on the line $x+y=t$ but, for $t\neq 0$, this level set is different of this line.

 So we have proved the following theorem.

\begin{thm} There are non-flat 3-webs which have a circular symmetry around some point.\end{thm}

\noindent {\bf Remarks}.

Using \cite{DJ}, we could obtain analogous results for local ($p+1$)-webs of codimension $n$ on a $p.n$ dimensional space.

We fall on this results  trying to construct a counter-example to the celebrate "Gronwall's conjecture for planar $3$-webs" first appearing in \cite{G} (see \cite{SA} for the better history of this conjecture and the last results). To be more precise, we tried to construct a {\bf linear} non-flat 3-web, with a symmetry around some point, which is not a homography. In fact, after very long calculations via Maple, it seems that there is no example of such 3-web.

\end{document}